\documentclass{article}
\usepackage{graphicx}

\newtheorem{theorem}{Theorem}
\newtheorem{lemma}[theorem]{Lemma}

\newtheorem{cor}[theorem]{Corollary}

\title{New principles for auxetic periodic design}
\author{Ciprian S. Borcea and Ileana Streinu}

\begin{document}
\maketitle

\begin{abstract}
We show that, for any given dimension $d\geq 2$, the range of distinct possible designs for
periodic frameworks with auxetic capabilities is infinite. We rely on a purely geometric approach to auxetic trajectories developed within our general theory of deformations of periodic frameworks. 
\end{abstract}



\section{Introduction}

New digital manufacturing techniques have vastly expanded the possibilities of generating complex three-dimensional structures, across length scales, and have opened up new opportunities for kinematic and geometric design to address functional desiderata. This paper is concerned with periodic structures and {\em metamaterials} with {\em auxetic capabilities}, a challenging and fast evolving topic  at the intersection of mathematics, mechanical design  and 
materials science \cite{elipe:lantada:auxeticGeometries:2012,lee:singer:thomas:microNanoMaterials:advMat:2012,
reisEtAl:designerMatter:2015,borcea:streinu:geometricAuxetics:RSPA:arxiv:2015}.  Our contribution is to derive new principles for auxetic design from the geometric theory of auxetic deformations recently introduced in \cite{borcea:streinu:geometricAuxetics:RSPA:arxiv:2015}.

\paragraph{Auxetic behavior.} When stretched, most materials will shrink laterally. {\em Auxetic behavior} is the rather counter-intuitive property exhibited by some materials that widen laterally upon stretching. In elasticity theory, such materials are said to have {\em negative Poisson's ratios} \cite{greaves:surveyPoissonRatios:resNotesRoyalSoc:2013}. The promise of various applications and increased interest in obtaining synthetic structures or metamaterials with this type of response to tensile loading has led to a sequence of studies, with emphasis on cellular and periodic structures \cite{lakes:negativePoisson:1987,evans:etAl:molecularNetwork:Nature:1991,
milton:compositePoisson:JMechPhysSol:1992,grima:alderson:evans:auxeticRotatingUnits:2005,
greaves:lakes:etAl:PoissonRatio:2011,mitschkeEtAl:auxetic:rspa:2013,
mitschkeEtAl:JSS:2016,
tanaka:bistiffness:rspa:2013,cabras:brun:auxetic:rspa:2014}. However, the repertory of auxetic designs proposed in the literature remained confined to a few dozen examples in dimension two and much less in dimension three \cite{elipe:lantada:auxeticGeometries:2012}. The authors of \cite{lee:singer:thomas:microNanoMaterials:advMat:2012} remark on p.4792 that \textit{``it has been a challenge to design 3D auxetic micro-/nano-structured materials''.}

\paragraph{New foundations for periodic auxetics.} In \cite{borcea:streinu:geometricAuxetics:RSPA:arxiv:2015}, we introduced a purely geometric approach to auxetic deformations for crystalline materials and man-made mechanical structures modeled as {\em periodic bar-and-joint frameworks.} This approach, reviewed in Section \ref{sec:auxeticTheory} below, presents a number of distinct advantages over the conventional route through Poisson's ratios. First of all, as a rigorous mathematical theory, the model can be applied to a wide range of structures, across length scales,  provided that a periodic bar-and-joint framework organization is the dominant feature. There is no need for experimental or simulated determinations of Poisson's ratios, since auxetic capabilities can be recognized by strictly geometric criteria. In fact, our mathematical theory works in arbitrary dimension. Moreover, the geometric approach clarifies the analysis of flexible structures with several degrees of freedom. In this case, certain deformation trajectories may be auxetic, while many other trajectories would not be auxetic. Thus, the notion of auxetic behavior must refer only to a certain type of one-parameter deformations of a periodic structure.
 
\begin{figure}[h]
\centering
 {\includegraphics[width=0.6\textwidth]{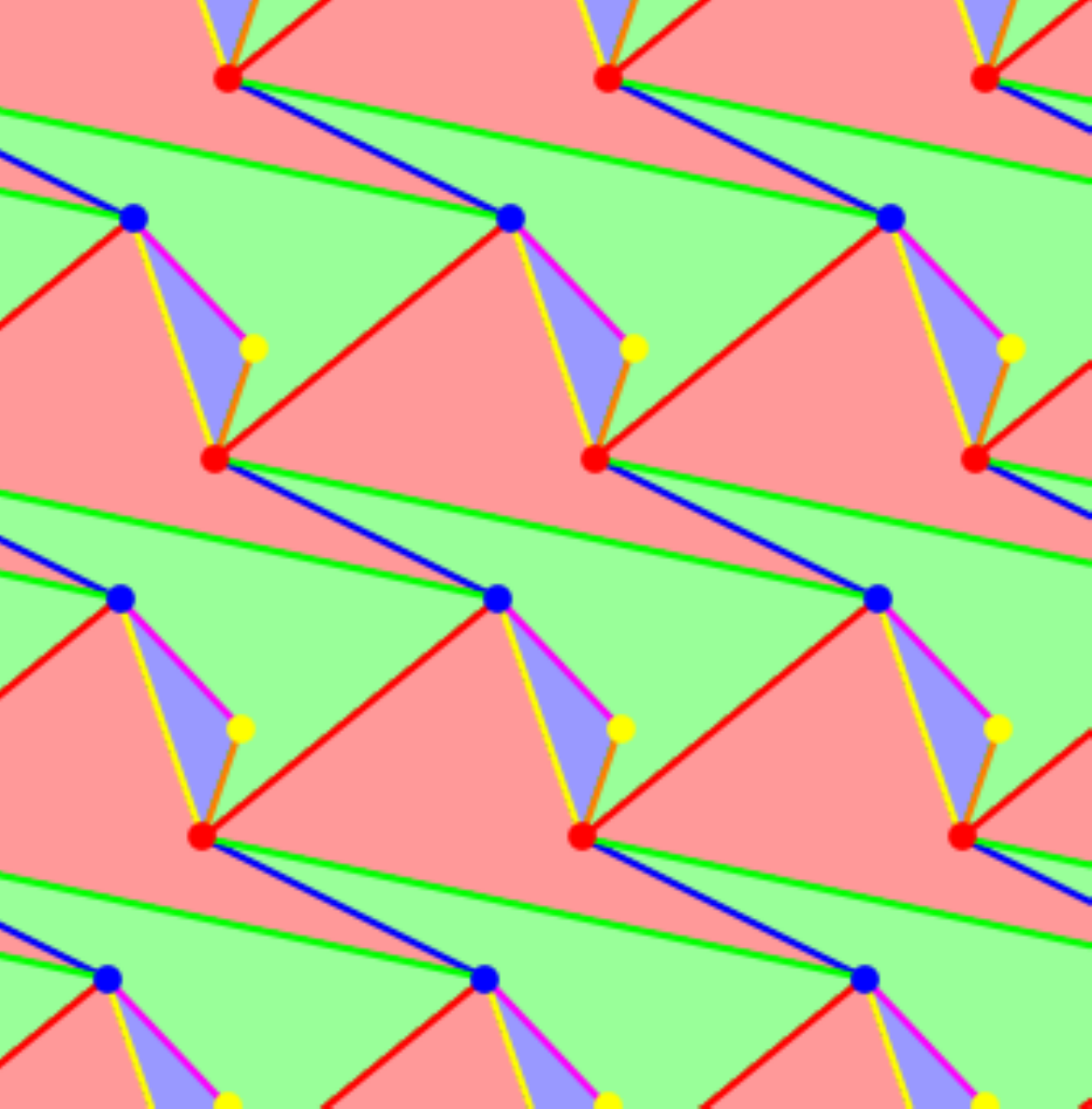}}
 \caption{A planar periodic pseudo-triangulation.}
 \label{FigPT}
\end{figure}

\paragraph{Auxetic design.} The ascendancy of the geometric approach is probably most conspicuous in matters of design. We have shown in \cite{borcea:streinu:kinematicsExpansive:ark14:2014,borcea:streinu:liftingsStresses:dcg:arxiv:2015} that the stronger notion of {\em expansive behavior}, when all distances between joints increase or stay the same, can be completely elucidated in dimension two in terms of a class of periodic structures called \textit{periodic pseudo-triangulations} (and kinematic equivalence classes of refinements to pseudo-triangulations). An example is presented in Figure~\ref{FigPT}.

\noindent
We have also shown that \textit{expansive implies auxetic,} hence this leads to an infinite gallery of planar auxetic periodic mechanisms, by virtue of the fact that all periodic pseudo-triangulations have exactly one degree of freedom to deform. While periodic pseudo-triangulations are easy to generate (see \cite{borcea:streinu:geometricAuxetics:RSPA:arxiv:2015} section 5(i) for a description and illustration of the procedure), proving 
their stated properties is not so elementary. Interested readers can find the full treatment in \cite{borcea:streinu:liftingsStresses:dcg:arxiv:2015}.

The {\em expansive implies auxetic principle} is valid in arbitrary dimension, but the  structure of expansive periodic frameworks in three or higher dimensions is only partially understood \cite{borcea:streinu:expansivePeriodic:IMAoxford:arxiv:2015}. In \cite{borcea:streinu:geometricAuxetics:RSPA:arxiv:2015}, we have relied on the suggestive value of necessary conditions for expansiveness for a couple of new designs of three dimensional periodic frameworks with auxetic capabilities.

\paragraph{Main contribution.} In the present work, we formulate and prove, in arbitrary dimension $d\geq 2$, a general principle for \textit{converting a finite linkage with adequate prerequisites on $d$ pairs of unconnected joints into a  periodic framework with auxetic capabilities.} From the standpoint of geometric auxetics \cite{borcea:streinu:geometricAuxetics:RSPA:arxiv:2015} these prerequisites are natural, elementary and easily satisfied.  This implies endless possibilities for auxetic design in arbitrary dimension.

\paragraph{Overview.} In order to control the number of degrees of freedom, we start with
a finite linkage in $R^d$ without self-stress, that is, with infinitesimally independent edge constraints. If the linkage has $d$ pairs of vertices which provide, as free vectors, a basis of $R^d$, we adopt the lattice generated by these vectors as {\em periodicity lattice} and obtain an
associated periodic framework by replicating the finite piece via all translations in the periodicity
lattice. This basic construction converts the finite linkage into a periodic framework with the
same number of (periodic) degrees of freedom, or, more precisely, with the same smooth local deformation space. Thus, in order to obtain auxetic capabilities for the periodic framework
we have to require adequate deformations for the finite linkage. It will be seen,
as we unfold this scenario, that the resulting requirement is simply an expression of the defining property of auxetic trajectories and can be
satisfied by an infinite gallery of finite linkage designs. We emphasize the fact that both the
simplicity of the principle and the unlimited variety of possible examples come from our
geometric theory of auxetic deformations \cite{borcea:streinu:geometricAuxetics:RSPA:arxiv:2015}.

\medskip
After reviewing in Section \ref{sec:auxeticTheory} the essential ingredients of this geometric approach, we describe the passage from finite linkages to periodic frameworks  in Section \ref{sec: F to P}. The design requirement for obtaining periodic frameworks with auxetic capabilities is stated in Section \ref{sec: principle} and then construction techniques for finite linkages with appropriate features are developed in Section \ref{sec: techne}. Proceedures for reducing the number of degrees of freedom to one and obtaining  auxetic periodic mechanisms are described next. In Section \ref{sec: affine} we show that affine transformations of a periodic framework preserve the infinitesimal auxetic cone. This result is relevant for classification criteria. A gallery of examples in Section \ref{sec: gallery} presents an infinite series of auxetic periodic mechanisms in dimension three. We mention some topics worthy of further investigation in our concluding Section \ref{sec: conclusion}.

\section{The geometric theory of periodic auxetics}
\label{sec:auxeticTheory}

In this section, we review the essential notions of our  geometric theory of auxetic one-parameter deformations of periodic frameworks. For more background and full details, we refer to \cite{borcea:streinu:periodicFlexibility:2010,borcea:streinu:geometricAuxetics:RSPA:arxiv:2015}.

\paragraph{Periodic graph.} A $d$-periodic graph is a pair $(G,\Gamma)$, where $G=(V,E)$ is a simple infinite graph with vertices $V$, edges $E$ and finite degree at every vertex, and $\Gamma \subset Aut(G)$ is a free Abelian group of automorphisms which has rank $d$, acts without fixed points on vertices and edges and has a finite number of vertex and orbits. The group  $\Gamma$ is isomorphic to $Z^d$ and is called the {\em periodicity group}  of the periodic graph $G$. Its elements $\gamma \in \Gamma \simeq Z^d$ are referred to as {\em periods} of $G$.

\paragraph{Periodic placement of a periodic graph.} A periodic placement of a $d$-periodic graph $(G,\Gamma)$ in $R^d$ is defined by two functions:

\[ p:V\rightarrow R^d \ \ \mbox{and} \ \ \pi: \Gamma \hookrightarrow {\textbf T}(R^d) \]

\noindent 
where $p$ assigns points in $R^d$ to the vertices $V$ of $G$ and $\pi$ is an injective homomorphism of $\Gamma$ into the group ${\textbf T}(R^d)$ of translations in the Euclidean space $R^d$, with $\pi(\Gamma)$ being a lattice of rank $d$. These two functions must satisfy the natural compatibility condition: $ p(\gamma v)=\pi(\gamma)(p(v)).$

\paragraph{Periodic framework.} A placement which does not allow the end-points of any edge to have the same image defines a {\em $d$-periodic bar-and-joint framework} in $R^d$, with edges $(u,v)\in E$ corresponding to bars (segments of fixed length) $[p(u),p(v)]$ and vertices corresponding to (spherical) joints. Two frameworks are considered equivalent when one is obtained from the other by a Euclidean isometry. 

\paragraph{Periodic deformation.} A {\em one-parameter deformation of the periodic framework} $(G,\Gamma, p,\pi)$ is a  (smooth) family of placements  $p_{\tau}: V\rightarrow R^d$  parametrized by time $\tau \in (-\epsilon, \epsilon)$ in a small neighborhood of the initial placement $p_0=p$, which satisfies two conditions: (a) it maintains the lengths of all the edges $e\in E$, and (b) it maintains periodicity under $\Gamma$, via faithful representations $\pi_{\tau}:\Gamma \rightarrow {\textbf T}(R^d)$ which {\em may change with $\tau$ and give an associated variation of the periodicity lattice of translations} $\pi_{\tau}(\Gamma)$. 

\medskip 
 After choosing an independent set of generators for the periodicity lattice $\Gamma$, the image $\pi_{\tau}(\Gamma)$ is completely described via the $d\times d$ matrix $\Lambda_{\tau}$ with column vectors given by the images of the generators under $\pi_{\tau}$. The associated {\em Gram matrix} is given by:

$$ \omega_{\tau}=\omega(\tau)=\Lambda^t_{\tau}\Lambda_{\tau}. $$

\paragraph{Auxetic path.} A deformation path $(G,\Gamma, p_{\tau},\pi_{\tau}), \tau\in (-\epsilon,\epsilon)$ is called auxetic when the curve of Gram matrices $\omega(\tau)$ defined above has all its
tangents in the cone of positive semidefinite symmetric $d\times d$ matrices. When all tangents are in the positive definite cone, the deformation is called {\em strictly auxetic}.

\medskip 
\noindent
This form of the auxeticity criterion, established in \cite{borcea:streinu:geometricAuxetics:RSPA:arxiv:2015}, Thm. 3.1, will be convenient for
our present purposes. We note that the auxetic character of a one-parameter deformation is  determined by the curve of Gram matrices of a basis of periods and strict auxeticity at one instance implies strict auxeticity in a neighborhood.

\section{From finite to periodic}\label{sec: F to P}

We assume a given dimension $d\geq 2$. A {\em linkage} in $R^d$ is a pair $L=(G,p)$, where
$G=(V,E)$ is a simple connected graph with $n=|V|$ vertices and $m=|E|$ edges and $p: V\rightarrow R^d$ is a placement of the vertices in $R^d$. Edges are then conceived as 
rigid straight bars between the corresponding vertices, which serve as spherical joints for the
linkage. It is assumed here that all edges correspond to non-zero segments. Linkages which differ by an isometry of $R^d$ are considered equivalent, that is $T\circ p$ and $p$ are equivalent placements for any Euclidean isometry $T$. Simply put, equivalent linkages are 
identified.

\medskip
The kinematics of linkages is a classical topic and we mention here only aspects and results directly relevant for our constructions. In particular, we use only {\em linkages with infinitesimally independent edge constraints}. In other words, the rows of the
 {\em rigidity matrix} are linearly independent. By the implicit function theorem, the local deformation space will be {\em smooth},
of dimension

\begin{equation}\label{Fdof}
f=nd -m- {d+1\choose 2}
\end{equation}

\noindent
where $ {d+1\choose 2}$ is the dimension of the Euclidean isometry group $E(d)$. We will need
$f \geq 1$, when $f$ is also called the number of {\em degrees of freedom} of the linkage.

\medskip
For the construction we are about to describe, we assume that $d$ pairs of vertices 
$\{ v_{i(k)},v_{j(k)} \} \subset V$, $k=1,...,d$, have been given, with the property that the vectors

\begin{equation}\label{basis}
\lambda_k=p(v_{j(k)})-p(v_{i(k)})\in R^d, \ \ \  k=1,...,d
\end{equation}

\noindent
form a basis of $R^d$. Note that, while in each pair $v_i\neq v_j$, different pairs may share
one vertex.

\medskip
Let $\Lambda=\{ \lambda=\sum_{k=1}^d n_k\lambda_k\ |\ n_k\in Z \}$, be the rank $d$ lattice generated
by this basis. This lattice will play the role of periodicity lattice for the periodic framework
associated to our linkage with marked pairs.

\medskip
We observe first that when we identify in $G$ all vertices which appear in one of the given pairs,
that is, if we put $v_{i{k}}\equiv v_{j(k)}$, $k=1,..,d$ and maintain all edges (some possibly converted into loops), we obtain a {\em quotient multigraph} with exactly $\tilde{n}=n-d$ vertices and $\tilde{m}=m$ edges.

\medskip
Let ${\cal V}\subset V$ be a complete set of representatives for the $n-d$ vertices of the quotient
multigraph. Then, the vertices $\tilde{V}$ of the periodic graph will be recorded as symbols
$ v+\lambda,\ \ v\in {\cal V}, \ \lambda \in \Lambda $
and placed by $\tilde{p}: \tilde{V} \rightarrow R^d$ at $\tilde{p}(v+\lambda)=p(v)+\lambda$.
The edge set $\tilde{E}$ of the periodic graph $\tilde{G}=(\tilde{V},\tilde{E})$ is made of all (formal) $\Lambda$ translates of $E$ (with $V$ seen as included in $\tilde{V}$). For 
periodicity group $\Gamma$ we have  $\Lambda$ itself, with the natural action on $\tilde{G}$.
The resulting periodic framework $(\tilde{G},\Gamma, \tilde{p},\pi)$, where $\pi$ is the identification  $\Gamma=\Lambda$, does not depend on the choice of representatives ${\cal V}$.

\begin{figure}[h]
\centering
 {\includegraphics[width=0.6\textwidth]{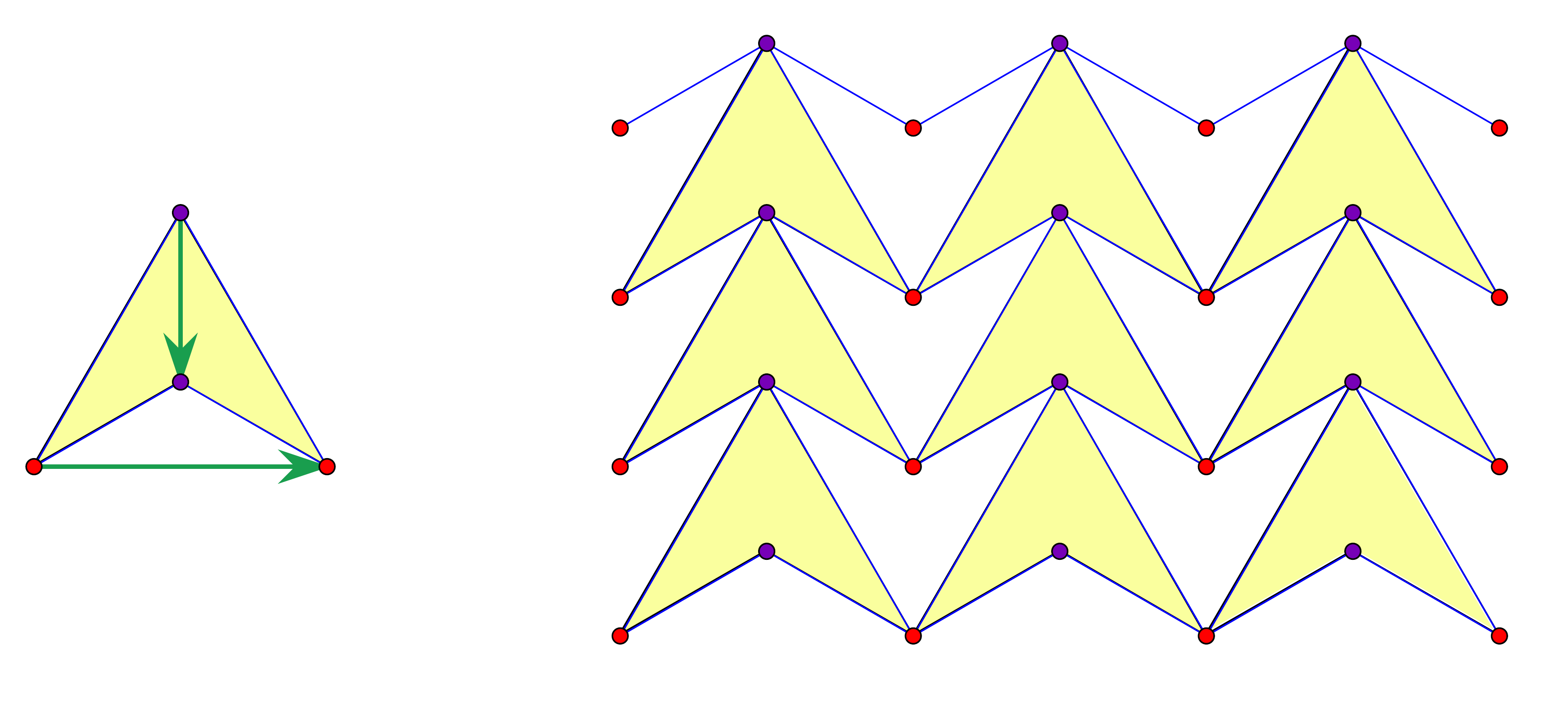}}
 \caption{A planar four-bar mechanism with a marked pair of vertices and the associated 
``double arrowhead" periodic framework.}
 \label{FigFtoP}
\end{figure}

Figure~\ref{FigFtoP} illustrates the procedure for $d=2$. The planar linkage is a four-bar mechanism configured as a pseudo-triangle. The chosen pairs of vertices correspond to the diagonals. The associated periodic framework is a periodic pseudo-triangulation described as a ``double arrowhead" pattern.

We summarize the general procedure as follows.

\begin{theorem}\label{fTOp}
Let $L=(G,p)$ be a connected linkage in $R^d$, with $n$ vertices and $m\leq nd-{d+1\choose2}$ edges, such that:

(a) the $m$ edges impose infinitesimally independent constraints,

(b) there are $d$ marked pairs of vertices, with the corresponding vectors (\ref{basis}) providing a basis of $R^d$.

\medskip
Then, there is a unique periodic framework $\tilde{L}=(\tilde{G},\Gamma,\tilde{p},\pi)$, with the following four properties:

\medskip \noindent
(i) the periodicity group
$\Gamma$ is identified with the lattice generated by the given basis, 

\noindent
(ii) the linkage $L$ is contained in the framework $\tilde{L}$,

\noindent
(iii) the edges of $L$ provide a complete set of representatives for the edge orbits of $\tilde{L}$ under periodicity,

\noindent
(iv) every vertex orbit of $\tilde{L}$  has at least one representative in $L$.

\medskip
This unique associated periodic framework $\tilde{L}$ has $\tilde{n}=n-d$ vertex orbits and $\tilde{m}=m$ edge orbits and the local deformation spaces of $L$ and $\tilde{L}$
can be identified. Their common dimension is

\begin{equation}\label{FPdof}
f=nd -m- {d+1\choose 2}=\tilde{n}d -\tilde{m}+{d\choose 2}
\end{equation}

\end{theorem}

\medskip 
The claim on preservation of degrees of freedom follows immediately from \cite{borcea:streinu:periodicFlexibility:2010}, p.2641.
The last term in (\ref{FPdof}) is the dimension formula for periodic deformations when 
edge constraints are infinitesimally (and hence locally) independent.

\medskip
The reader may observe that the planar framework in Figure~\ref{FigPT} can also be obtained
(in two ways) by a passage from finite to periodic.

\section{The main auxetic design principle}\label{sec: principle}

The correspondence described in Theorem~\ref{fTOp} turns a finite linkage $L=(G,p)$ into a 
periodic framework $\tilde{L}=(\tilde{G},\Gamma,\tilde{p},\pi)$ with the same local deformation space. Thus, a one-parameter local deformation $(G,p_{\tau})$ for $L$ turns into a 
(periodic) one-parameter local deformation $(\tilde{G},\Gamma,\tilde{p}_{\tau},\pi_{\tau})$ for
$\tilde{L}$. If we want the latter to be an {\em auxetic path}, the criterion of 
\cite{borcea:streinu:geometricAuxetics:RSPA:arxiv:2015}, Thm. 3.1, recalled above at the end of Section~\ref{sec:auxeticTheory}, requires the generators of the periodicity lattice to give a 
curve of Gram matrices $\omega(\tau)$ with tangent directions $\frac{d\omega}{d\tau} (\tau)$
{\em in the positive semidefinite cone}. Obviously, this curve is determined by the effect of the
linkage deformation $p_{\tau}$ on the $d$ pairs of vertices marked on the linkage at the outset.

\medskip
The essence of our auxetic design principle can be stated already.

\begin{theorem}\label{principle}
The periodic framework $\tilde{L}$ has a non-trivial auxetic deformation path if and only if the finite linkage $L$ has
a local one-parameter deformation with the following property:   the Gram matrix of the
basis given by the $d$ marked pairs of vertices evolves (in the space of $d\times d$ symmetric matrices) as a non-constant curve with all its tangents in the 
positive semidefinite cone.
\end{theorem}

\medskip
An important case, which involves only infinitesimal considerations, is when the tangent direction $\frac{d\omega}{d\tau} (0)$ at the initial moment $\tau=0$ is actually in the {\em positive definite cone}. By continuity, this is enough to guarantee an interval 
$\tau \in (-\epsilon, \epsilon)$ where the tangent remains in the positive definite cone and
the periodic deformation is strictly auxetic.

\medskip
For explicitness, we unfold the more formal details. 

\medskip
With notations introduced above in (\ref{basis}), the variation with $\tau$ of the marked basis
is given by

\begin{equation}\label{basisVar}
\lambda_k(\tau)=p_{\tau}(v_{j(k)})-p_{\tau}(v_{i(k)})\in R^d, \ \ \  k=1,...,d
\end{equation}

\noindent
The $d\times d$ matrix with these column vectors is denoted $\Lambda(\tau)$, hence the Gram matrix of the marked basis is

\begin{equation}\label{GramVar}
\omega(\tau)=\Lambda^t(\tau) \Lambda(\tau) .
\end{equation}

\noindent
The velocity vector at moment $\tau$ for this parametrized curve is

\begin{equation}\label{velocity}
\frac{d\omega}{d\tau}(\tau)=\frac{d\Lambda^t}{d\tau}(\tau) \Lambda(\tau)+ 
\Lambda^t(\tau) \frac{d\Lambda}{d\tau}(\tau) .
\end{equation}

\medskip \noindent
The auxeticity requirement is that all velocity vectors (\ref{velocity}) for $\tau$ in a small neighborhood of 0 belong to the positive semidefinite cone. At moment $\tau=0$, we have
$\Lambda(0)=\Lambda$ and

\begin{equation}\label{velocity(0)}
\frac{d\omega}{d\tau}(0)=\frac{d\Lambda^t}{d\tau}(0) \Lambda+ 
\Lambda^t \frac{d\Lambda}{d\tau}(0) .
\end{equation}

\medskip \noindent
with the derivative $\frac{d\Lambda}{d\tau}(0)$ depending only on the infinitesimal deformation
corresponding to $p_{\tau}$. As noticed above, when (\ref{velocity(0)}) is in the interior of the
positive semidefinite cone, that is $\frac{d\omega}{d\tau}(0)$ is positive definite, strict auxeticity
follows for small enough $\tau$. We state explicitly this corollary, since this form of the principle implicates only infinitesimal deformations of $L$ and is most useful for constructing
examples.

\begin{cor}\label{strict}
If, for some infinitesimal deformation of $L$, the corresponding infinitesimal variation
$\frac{d\omega}{d\tau}(0)$ of the Gram matrix is positive definite, the periodic framework $\tilde{L}$ has a strictly auxetic deformation path.
\end{cor}
 
\begin{figure}[h]
\centering
 {\includegraphics[width=0.5\textwidth]{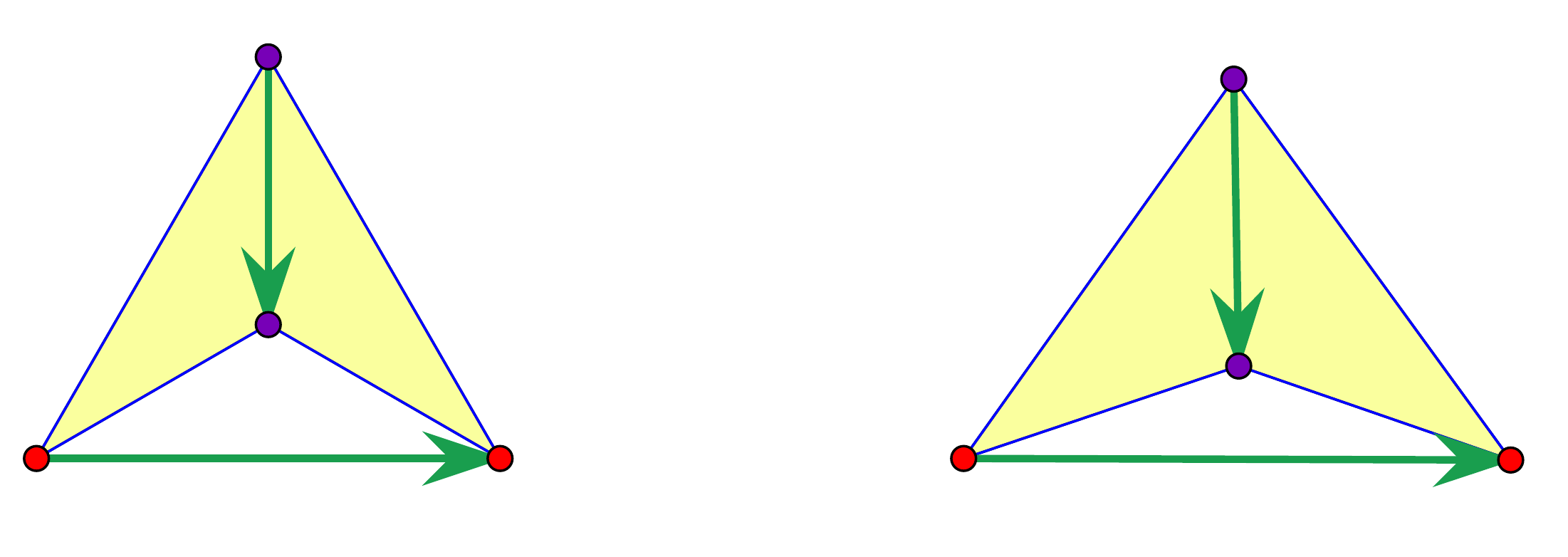}}
 \caption{Deforming the quadrilateral in Figure~\ref{FigFtoP}. }
 \label{FigQuad}
\end{figure}

For a simple illustration, we revisit the example given in Figure~\ref{FigFtoP}. The quadrilateral
has one degree of freedom and the obvious local deformation for our symmetric configuration,
shown in Figure~\ref{FigQuad}, increases the lengths of the two diagonals, while maintaining their orthogonality. This means that $\frac{d\omega}{d\tau}(0)$ is a diagonal $2\times 2$ matrix
with positive diagonal entries. Hence, the double arrowhead periodic framework is a strictly auxetic one degree of freedom periodic mechanism (as long as the quadrilateral remains concave
i.e. remains a pseudo-triangle).

\section{Construction techniques}\label{sec: techne}

In this section we discuss construction techniques for obtaining finite linkages $L$
which satisfy the auxetic requirement formulated in Theorem~\ref{principle}, or rather the strict
auxetic criterion of Corollary~\ref{strict}. It will soon become apparent that examples can be
constructed {\em ad libitum} in any dimension $d\geq 2$. The general construction ideas surveyed here offer wide opportunities for applications, since additional specifications can be met by ingenuity and refinement in the finite linkage design.

\medskip
We begin with an examination of the case when all $d$ pairs of vertices marked on the linkage $L=(G,p)$ have a common vertex $v_0$. Thus, convenient labeling will have our basis expressed as

\begin{equation}\label{simplex}
\lambda_k=p(v_k)-p(v_0), \ \ k=1,...,d
\end{equation}

\medskip \noindent
with the Gram matrix $\omega=( \langle \lambda_i, \lambda_j \rangle )_{ij}$, $1\leq i,j \leq d$.
Thus, we aim at designing $L$, so that it has an infinitesimal deformation making the corresponding infinitesimal variation $\dot{\omega}=\frac{d\omega}{d\tau}(0)$ of the Gram matrix positive definite.

\medskip
We may assume, without violating the condition on infinitesimally independent edge constraints, that $L$ contains a rigid part (e.g. a $d$-dimensional simplex) and we shall refer to it as the
``scaffold". Then, we can attach to this scaffold, without any redundancy of new bars, other
elements of the linkage, with adequate control on the possible motion of certain ``vertices of interest". In our case, the vertices of interest are those labeled $v_0,...,v_d$, and we may connect them to the scaffolding as follows: $v_0$ is rigidly connected, while each $v_k$, $k=1,...,d$
is attached via a ``hinge connection", to be described in the next lemma, which allows only one degree of freedom relative to the scaffold, namely the rotation of $v_k$ around that hinge. 

\begin{lemma}\label{hinge}
Suppose we have a (spanning and minimally) rigid linkage in $R^d$, referred to as a ``scaffold".
Suppose we have a new vertex $v$ and we want to connect it to the scaffold in such a way that
it has only one  degree of freedom with prescribed direction for its infinitesimal motion relative to the scaffold.
Then, we choose a $(d-2)$-simplex, referred to as a ``hinge" and  position it in the hyperplane through $v$ with normal direction prescribed  by the infinitesimal motion allowed for $v$ (but away from $v$). We connect this hinge rigidly (and without redundancy) to the scaffold and connect $v$ to the $(d-1)$ vertices of the simplex. 
The resulting linkage allows only a rotation of $v$ around the hinge as relative motion and the infinitesimal direction of motion is the one prescribed.
\end{lemma}

\begin{figure}[h]
\centering
 {\includegraphics[width=0.4\textwidth]{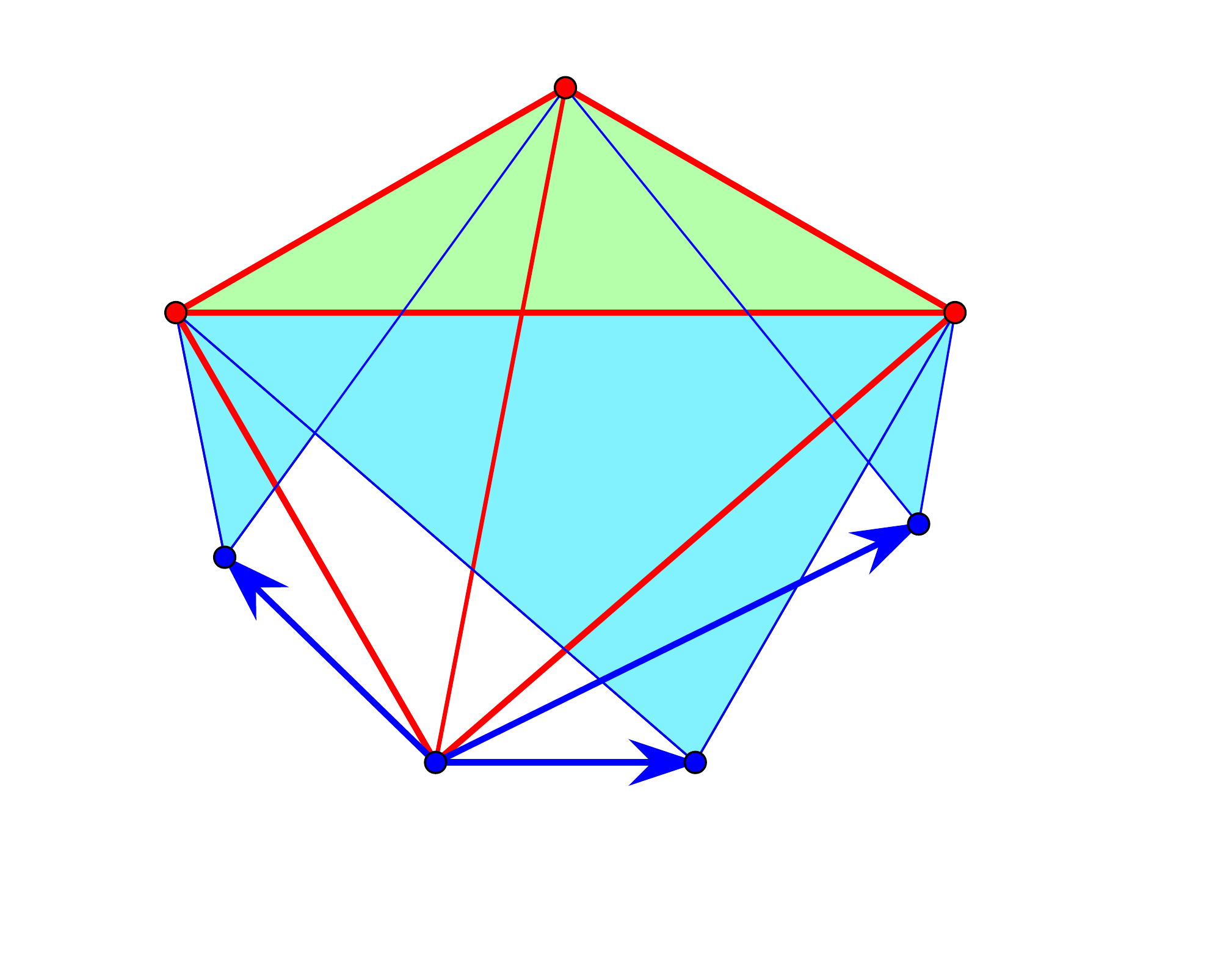}}
 \caption{Paneled simplex}
 \label{FigPS}
\end{figure}

\medskip \noindent
{\bf Remark.}\ In more suggestive terms, this lemma may be called ``the trapdoor principle",
since $v$ and the hinge form a panel (codimension one ``trapdoor") and this panel can only rotate around its fixed hinge. In Figure~\ref{FigPS} we show the result of applying this principle to
obtain a linkage in $R^3$ with controlled motion for three vertices relative to a (fixed) tetrahedral
scaffold. The three marked arrows would be $\lambda_k$, $k=1,2,3$ and this shows that 
their infinitesimal variation can be arranged to yield a positive definite $\dot{\omega}$, as argued 
in the next lemma.

\begin{lemma}\label{altitudes}
Suppose $p_0,p_1,...,p_d$ are points in $R^d$ with

$$ \lambda_k=p_k-p_0, \ \ k=1,...,d $$

\noindent
a basis of $R^d$. The point $p_0$ is assumed fixed and may be taken as the origin.
The points $p_k$ are subject to infinitesimal displacements $\mu_k\neq 0$ which have the direction of the corresponding altitude from $p_k$ in the simplex $[p_0...p_d]$ and are all pointing outwards. Then, the resulting infinitesimal variation $\dot{\omega}$ of the Gram matrix $\omega=( \langle \lambda_i, \lambda_j \rangle )_{ij}$ is positive definite.
\end{lemma}

\medskip \noindent
{\em Proof.}\ Let us assume that only one point, say $p_k$, moves infinitesimally by $\mu_k\neq 0$,
 with the other points fixed. Then $\dot{\omega}$ is positive semidefinite of rank one, with the only non-zero entry at $(k,k)$. The lemma follows by linear combination (with all coefficients one).

\begin{figure}[h]
\centering
 {\includegraphics[width=0.6\textwidth]{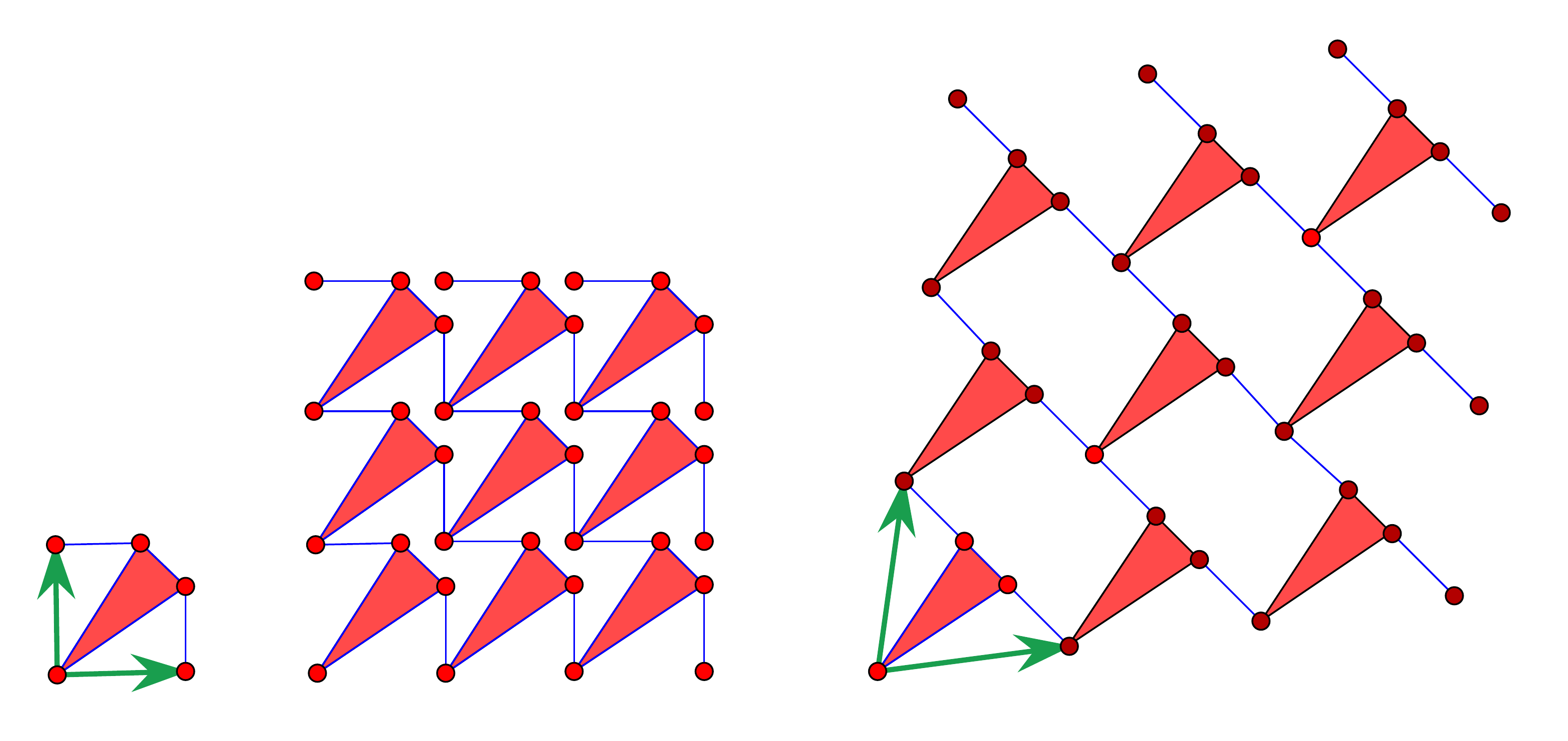}}
 \caption{A finite linkage with two degrees of freedom with auxetic capabilities for the associated periodic framework. The deployed configuration on the right can be reached via an auxetic path.}
 \label{Fig2deploy}
\end{figure}

\medskip
In summary, the trapdoor principle (Lemma~\ref{hinge}) shows how to design a  linkage
with prescribed infinitesimal motions for marked vertices $v_0,...v_d$. When these prescriptions are as in Lemma~\ref{altitudes} (with $p_i=p(v_i)$) we obtain a positive definite
infinitesimal variation of the Gram matrix of periods, hence strictly auxetic capabilities for
the associated periodic framework. Figure~\ref{Fig2deploy} illustrates this type of construction for $d=2$. 

\begin{figure}[h]
\centering
 {\includegraphics[width=0.8\textwidth]{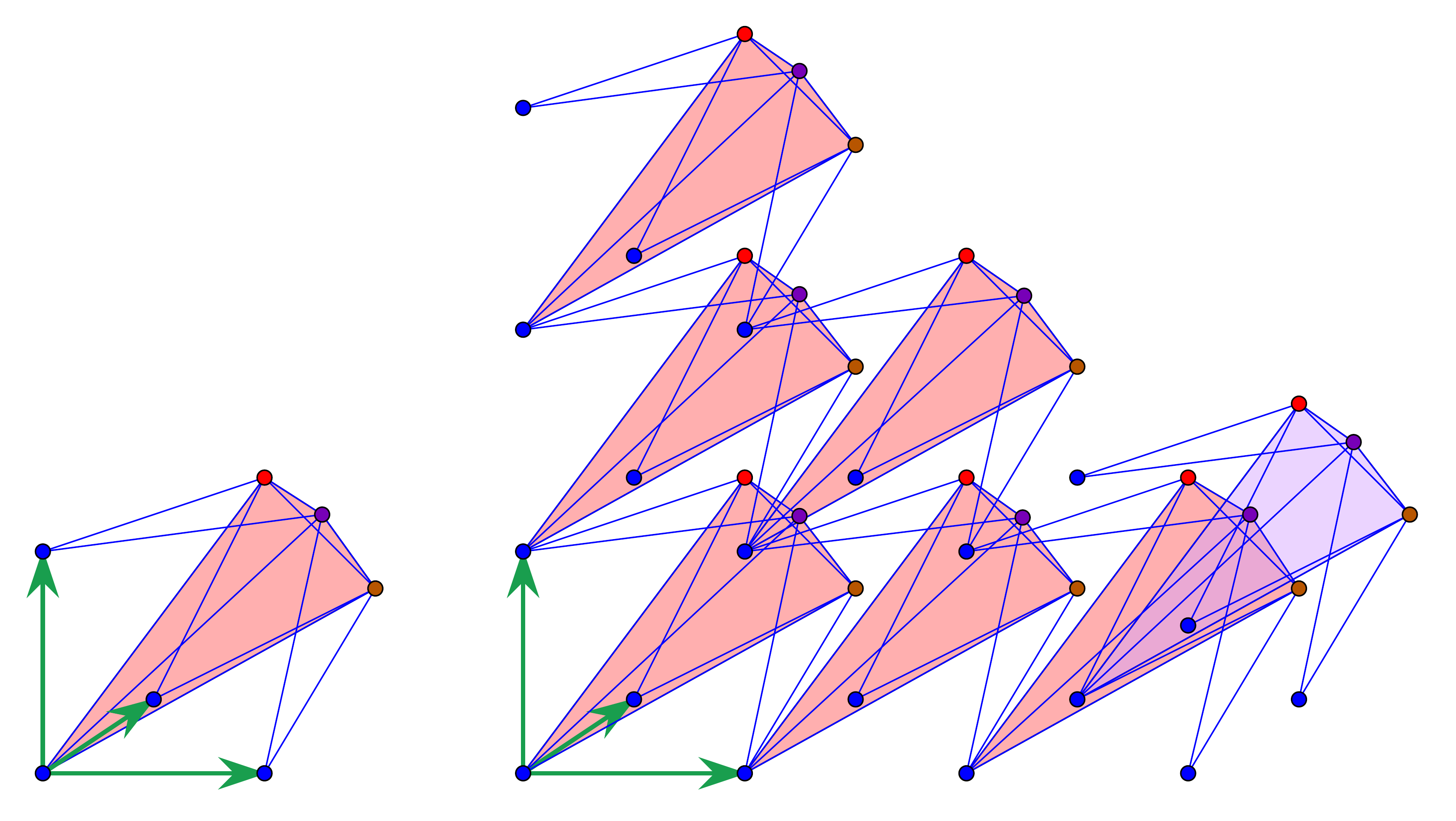}}
 \caption{The $d=3$ version of the planar scenario in Figure~\ref{Fig2deploy}. The finite linkage is a paneled tetrahedron, as in Figure~\ref{FigPS}. Only one `in depth' translate is shown.}
 \label{Fig3deploy}
\end{figure}

\medskip
Figure~\ref{Fig3deploy} shows the three-dimensional version of the same scheme.

Note that these constructions lead to linkages with $d$ degrees of freedom. In Section~\ref{sec: reduction} we describe a simple reduction procedure to one degree of freedom and strictly auxetic motion for the associated periodic framework. 

\medskip
It is fairly transparent by now that similar scenarios apply for other patterns of $d$ pairs of vertices in the finite linkage $L$. If we consider, for instance, the case of no common vertex
for any two pairs, we may label the pairs $(v_k,w_k)$, with basis 

$$ \lambda_k=p(w_k)-p(v_k), \ \ k=1,...,d. $$

\medskip
We may design $L$ with all $v_k$, $k=1,...,d$ fixed to the scaffold and $\lambda_k$ mutually
orthogonal and each $w_k$, constrained by hinge connections to move infinitesimally along
$\lambda_k$. Again, we obtain associated periodic frameworks with strictly auxetic capabilities.

\medskip
For other patterns of $d$ marked pairs one may use orthogonal splittings of $R^d$ and
maintain orthogonality for adequate partitions of the $d$ pairs. We dispense with further inventory pursuits here and consider instead a type of structure which will be useful for
deriving an explicit infinite series of non-isomorphic $d$-periodic graphs underlying auxetic periodic mechanisms for $d\geq 3$.

\begin{figure}[h]
\centering
 {\includegraphics[width=0.8\textwidth]{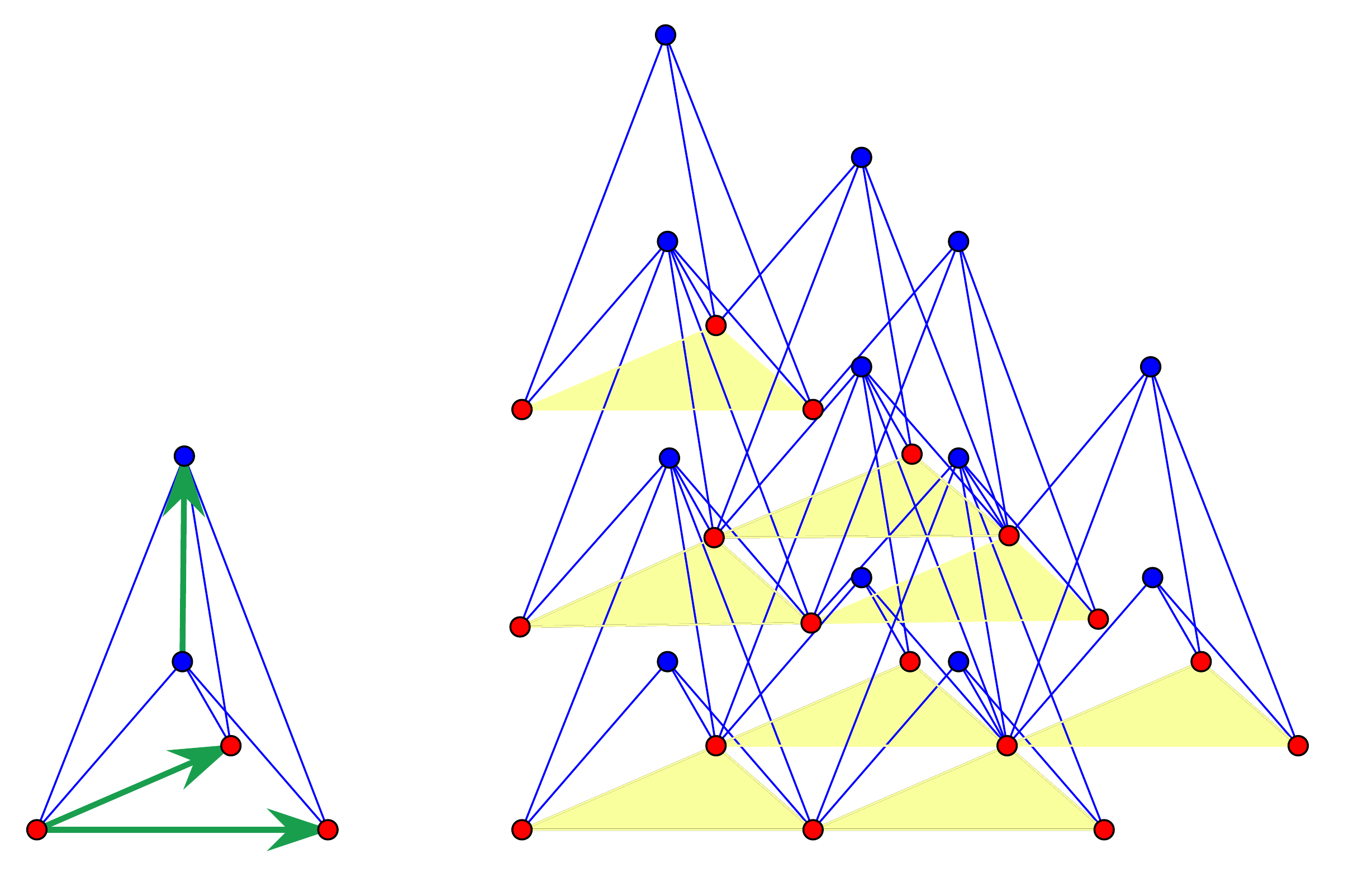}}
 \caption{A finite linkage  in $R^3$, with 5 vertices, 6 edges and 3 pairs of vertices marked by arrows. A fragment of the associated periodic framework is shown nearby. It has 2 orbits of vertices and 6 orbits of edges modulo periodicity.}
 \label{FigCadelniza}
\end{figure}

\medskip
Let $d\geq 3$ and consider the finite linkage in $R^d$ with $n=d+2$ vertices and $m=2d$ edges which corresponds to the case shown in Figure~\ref{FigCadelniza} for $d=3$. For a more suggestive description, we adopt the following language. There is a hyperplane, to be called `floor', which contains $d$ of the placed vertices, configured as a regular $(d-1)$-simplex, but not connected by edges. The two remaining vertices lie above the center of the $(d-1)$-simplex and both are connected by edges to all $d$ points on the floor. As a result, all bars from one or the other of the two points are of equal length and the line through the two points is orthogonal to the floor (and will be called `vertical'). The vectors given by $d$ marked pairs of vertices are: the vertical vector between the two points above the floor and the $(d-1)$ `horizontal' vectors from one point on the floor to the remaining $(d-1)$ points on the floor.

\noindent
The indicated vectors become periods and generate the periodicity lattice for the associated
periodic framework. We see that, when using horizontal periods, the floor goes to itself and
gives a floor for the periodic framework. Repeated translations by the vertical period will generate
an infinite array of such floors. Floors do not contain edges. In Figure~\ref{FigCadelniza}, replicas
of the floor $(d-1)$-simplex are highlighted. 

\medskip
There are ${d\choose 2}$ degrees of freedom and a natural set of parameters for the deformations of the finite linkage would be the (squared) distances between the vertices of the
floor simplex. A {\em strictly auxetic deformation trajectory} for the periodic framework can be immediately proposed based on a dilation motion for the floor simplex. 
The vertical vector remains orthogonal to the floor since the two end-points project to the circumcenter of the floor simplex. The 
end-point closer to the floor approaches faster than the remote one, resulting in a (squared) length increase for the vertical period vector. Thus, $\dot{\omega}$ is clearly positive definite for
this kind of local deformation trajectory. Other auxetic trajectories (with $\dot{\omega}$ semipositive definite of rank two) will be described in the next section, in connection with
ways of reducing from  ${d\choose 2}$ degrees of freedom to one.

\section{Reduction to one degree of freedom}\label{sec: reduction}

In this section we show first that the periodic frameworks with $d$ degrees of freedom obtained
when using our auxetic design principle in the manner described in Section~\ref{sec: techne} can be turned into periodic frameworks with a single degree of freedom which are locally strictly
auxetic. 

\medskip  
More precisely, let us assume that we have constructed a finite linkage $L=(G,p)$ in $R^d$, with $d$ pairs of vertices marked as $(v_0,v_k)$, $k=1,...,d$. We assume at the same time, that the rigid
part of $L$ called scaffold is fixed and includes $v_0$, which is placed at the origin. The linkage
has $d$ degrees of freedom due to the hinge connections of the $d$ vertices
$v_k$, $k=1,...,d$, to the scaffold. 

\medskip
We let $p_k=p(v_k)$, $k=1,...,d$ denote the initial positions and consider a deformation path
with infinitesimal displacements $\dot{p}_k$ for the vertices $v_k$, as needed when satisfying a strictly auxetic prerequisite.  In particular, this is the setting illustrated above in connection with Figures \ref{Fig2deploy} and \ref{Fig3deploy}.

\begin{lemma}\label{reduction}
The linkage $L$ can be turned into a linkage $L^*$ which has a single degree of freedom
and a deformation path with the same infinitesimal displacements $\dot{p}_k$ for
the vertices $v_k$, $k=1,...,d$.
\end{lemma}

\medskip \noindent
{\em Proof.}\ We introduce a new vertex $w$ with a sufficiently general position $q$. We connect
$w$ with all $v_k$, $k=1,...,d$. Then, the infinitesimal displacement  $\dot{q}$ of $w$ which
is compatible with the infinitesimal displacements $\dot{p}_k$ is uniquely determined by the
linear system:

\begin{equation}\label{unique}
\langle \dot{p}_k-\dot{q}, p_k-q \rangle=0, \ \ k=1,...,d
\end{equation}

\medskip \noindent
Now, we may apply the `trapdoor principle' (Lemma~\ref{hinge}) and construct a hinge connection for $w$ to the scaffold, with $\dot{q}$ as the allowed direction of infinitesimal
displacement. This yields the desired linkage $L^*$ with one degree of freedom.

\begin{figure}[h]
\centering
 {\includegraphics[width=0.7\textwidth]{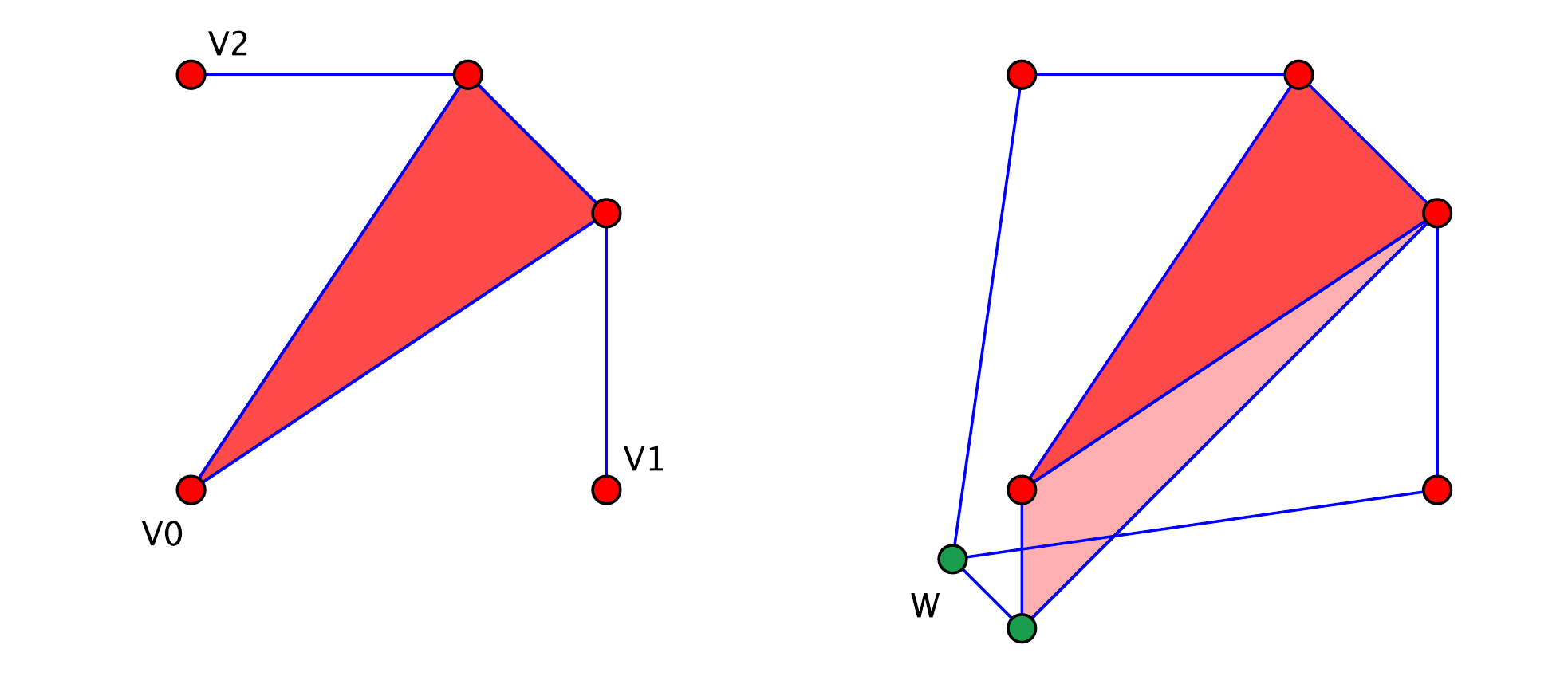}}
 \caption{The planar linkage in Figure~\ref{Fig2deploy} has two degrees of freedom but can be
converted to a single degree of freedom mechanism which retains the desired infinitesimal deformation on the vertices $v_0,v_1,v_2$.}
 \label{Fig2adapt}
\end{figure}
 
\medskip
We illustrate in Figure~\ref{Fig2adapt} a simple conversion of the linkage used earlier in Figure~\ref{Fig2deploy} into a linkage with just one degree of freedom. Obviously, this kind of 
operation has considerable leeway.  

\medskip \noindent
{\bf Remark.}\ For the associated periodic framework we obtain strict auxeticity (i.e. $\dot{\omega}$ positive definite) at the initial moment, hence a strictly auxetic deformation path defined on some interval $(-\epsilon,\epsilon)$.
 
\medskip
A somewhat different procedure will be described presently for the type of frameworks related to
Figure~\ref{FigCadelniza}. In this setting, $d\geq 3$ and the task is to reduce the degrees of freedom from ${d\choose 2}$ down to one. For simplicity and the benefit of figures, we conduct
our discussion in {\em dimension three}. The arguments in higher dimensions are completely analogous.   

\begin{figure}[h]
\centering
 {\includegraphics[width=0.9\textwidth]{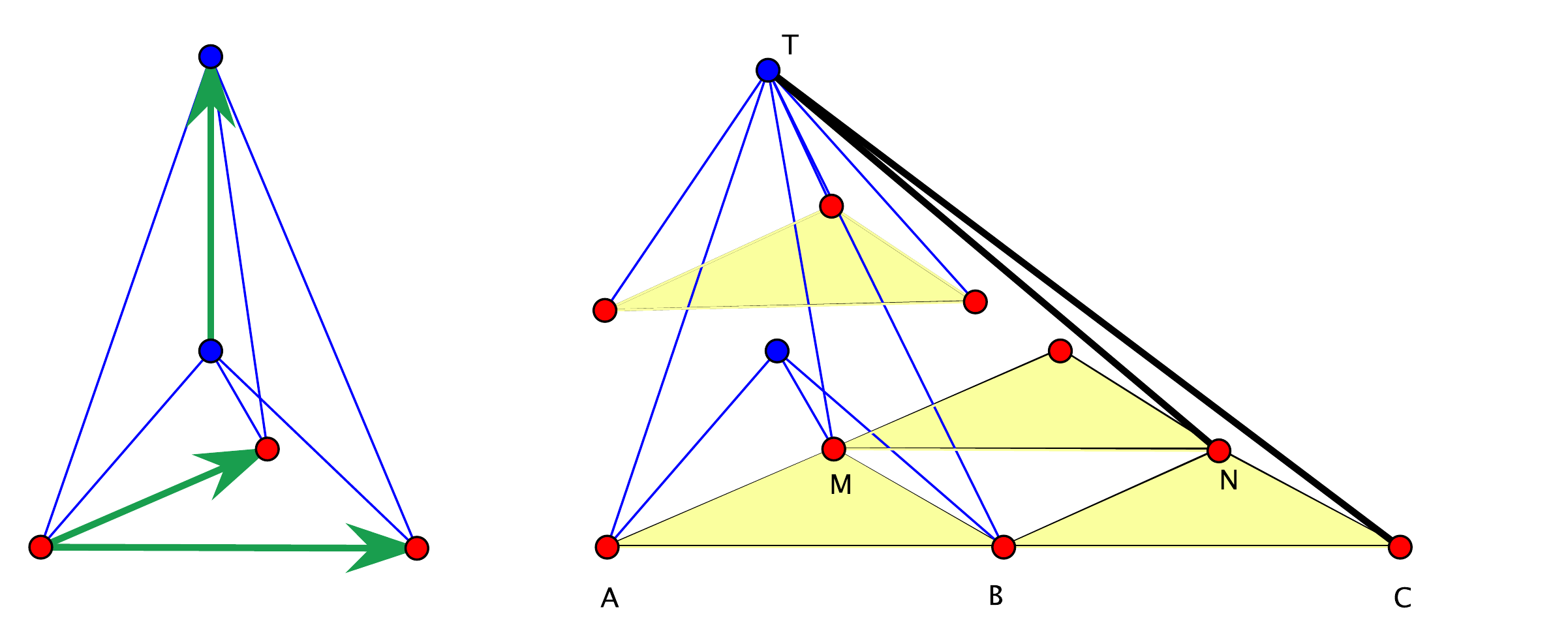}}
 \caption{Two new edge orbits in the associated periodic framework.}
 \label{Fig:roofing}
\end{figure}

\medskip 
The main idea is illustrated in Figure~\ref{Fig:roofing}. We operate in the associated periodic framework, where we want to introduce ${d\choose 2}-1$ new edge orbits, that is (for $d=3$)
two new edge orbits. Recall that highlighted triangles belong to stacked floors. Floors do not contain edges, but are organized by horizontal periodicity. The figure shows three triangles in
one floor and one triangle in the floor above. The representatives $TN$ and $TC$ of the two new edge orbits are shown as darker segments.

\medskip \noindent
Note that the three edges from the top vertex to the triangle in the upper floor belong, by periodicity, to the original periodic framework, which has six edges incident to every vertex.
Recall that we assumed an initial configuration with equilateral floor triangles and the vertical
periods positioned over the centers of floor triangles.
With the two new edge orbits, we obtain a periodic framework with the same periodicity lattice,
but with eight edges incident to every vertex. In the figure, only the top vertex shows all
eight bars incident to it. With these caveats taken into account, we proceed with the arguments
which prove the auxetic character of the resulting periodic mechanism. 

\medskip 
The triangle $TAC$ is determined by the two edges $TA$ and $TC$, together with the median
$TB$. This means that the horizontal period $AB$ has determined length. Since $MN$ is the same
period, the triangle $TMN$ is completely determined. Thus, the periodic mechanism can only
open up the dihedral angle of the planes $TAB$ and $TMN$ which have a common line with the direction of the period $AB$. Since $T$ remains on the perpendicular bisector of $AB$, so does $M$. Thus, $N$ remains on the perpendicular bisector of $BC$. We have argued earlier that the
vertical period remains vertical (i.e. perpendicular on the floor) and increases in length. All in all,
the infinitesimal variation of the Gram matrix of periods has two positive diagonal entries
(corresponding to the vertical period and $AM$) and zero elsewhere. Hence the motion is auxetic.

\begin{figure}[h]
\centering
 {\includegraphics[width=0.65\textwidth]{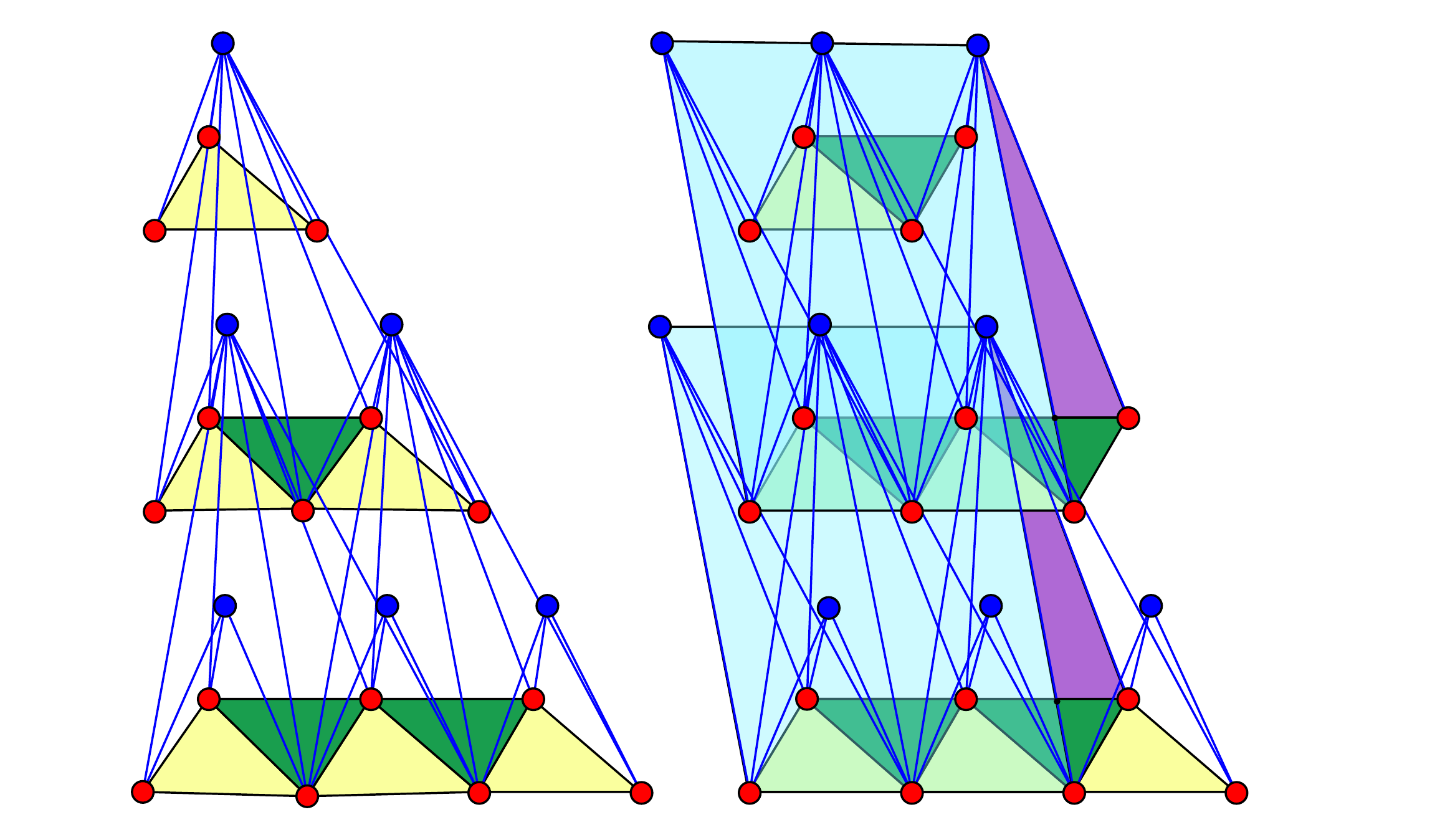}}
 \caption{Auxetic mechanism. The fragment shows three floors. Alternative view as  ``breathing stacked roofs". {\em Inhaling is auxetic.} The one degree of freedom deformation can be parametrized by the dihedral angle of a roof. The floors vary accordingly. The framework edges are all contained
in roof planes. Floors contain no edges but must respect periodicity constraints.}
 \label{FigMechanism}
\end{figure}

\medskip \noindent
{\bf Remark.}\ For a suggestive reading of Figure~\ref{Fig:roofing} in higher dimension $d$, the
$AB$ part of the floor triangle should be conceived as a facet of the floor simplex. The vertex $T$
will be connected to all but one of the vertices and edge-midpoints of the duplicated simplex
in the lower floor. All edges incident to $T$ in the resulting one degree of freedom periodic
framework are contained in two hyperplanes (`roofs'), which play the role of the two planes $TAC$ and $TMN$ in the figure. The auxetic mechanism may be fancied as a periodic array of
(crisscrossing and interlaced) opening roofs as proposed in Figure~\ref{FigMechanism}.

\begin{figure}[h]
\centering
 {\includegraphics[width=1\textwidth]{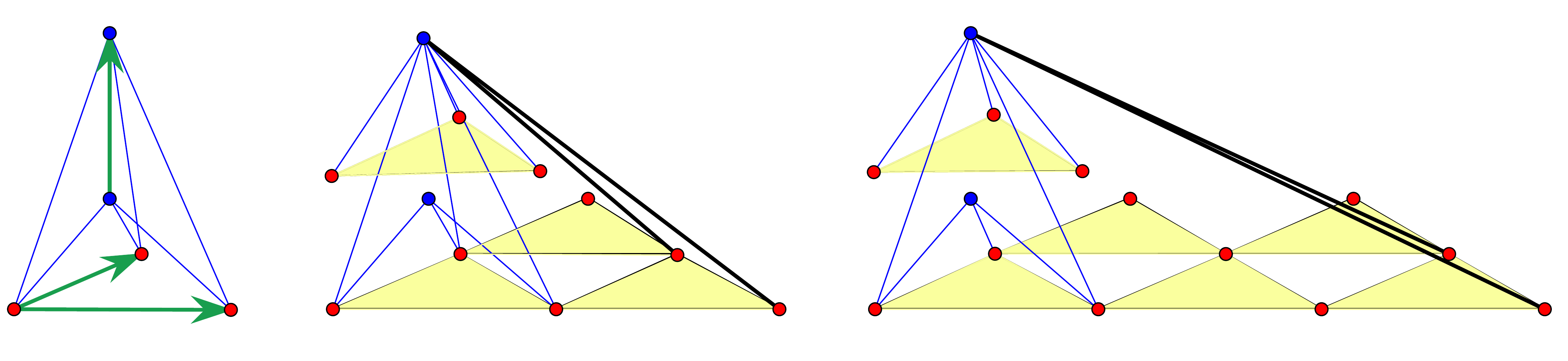}}
 \caption{Alternative way for introducing new edge orbits (in the same roof planes).}
 \label{FigSeries}
\end{figure}

\medskip
We may observe at this point that there are other ways of introducing new edge orbits for
a reduction to one degree of freedom mechanisms. If we follow the pattern indicated in Figure~\ref{FigSeries}, we obtain a series of non-isomorphic periodic graphs with 
framework realizations which are auxetic periodic mechanisms. Since the roof planes are the same in this series, the resulting auxetic mechanisms are kinematically equivalent, although
structurally distinct. This shows that equivalence criteria in periodic auxetics must be introduced
with some care. Another aspect deserving attention when attempting classifications is affine equivalence, discussed in the next section.

\section{Affine transformations}\label{sec: affine}

In this section we show  that the  infinitesimal auxetic cone of a periodic framework is preserved under affine transformations, that is, the natural isomorphism of the corresponding vector spaces of infinitesimal deformations \cite{borcea:streinu:periodicFlexibility:2010}, Prop. 3.7, pg. 2639, takes one auxetic cone to the other. While this fact is not directly intuitive,  it has a 
straightforward and short proof in our geometric theory of periodic auxetics.
 
\medskip
We adopt the following notations and setting: ${\cal F}=(G, \Gamma,p,\pi)$ is a periodic framework in $R^d$, with $n$ vertex orbits and  $m$ edge orbits. After a choice of independent generators for $\Gamma$, the lattice of periods $\pi(\Gamma)$ of the framework 
is described by
a $d\times d$ matrix $\Lambda$, with columns given by the images of the generators.
We denote by $p_i$, $i=1,...,n$ the positions of a complete set of representatives for the vertex orbits. An infinitesimal deformation of the framework  ${\cal F}$ is determined by the infinitesimal variation $\dot{p}_i$ of the positions $p_i$ and the infinitesimal variation $\dot{\Lambda}$
of the periodicity matrix $\Lambda$.

\medskip
Since the effect of translations is trivial, we assume our affine transformation to be a linear map $A: R^d\rightarrow R^d$. Then the transformed framework $A({\cal F})$ has vertex representatives at $Ap_i$ and periodicity matrix $\Lambda_1=A\Lambda$. The natural isomorphism between infinitesimal deformations takes $(\dot{p}_i,\dot{\Lambda})$
for ${\cal F}$ to $((A^t)^{-1}\dot{p}_i, (A^t)^{-1}\dot{\Lambda})$ for $A({\cal F})$. Thus, $\dot{\Lambda}_1=(A^t)^{-1}\dot{\Lambda}$ and 
$\dot{\Lambda}_1^t= \dot{\Lambda}^t A^{-1}$.

\medskip
Auxetic deformations for ${\cal F}$ are those with positive semidefinite $\dot{\Lambda}^t\Lambda+\Lambda^t \dot{\Lambda}$ and the isomorphism gives

\begin{equation}\label{invariance}
  \dot{\Lambda}_1^t\Lambda_1+\Lambda_1^t \dot{\Lambda}_1
=\dot{\Lambda}^t\Lambda+\Lambda^t \dot{\Lambda} \succeq 0 
\end{equation}

\noindent
confirming the preservation of the auxetic cone.

\medskip \noindent
{\bf Remark.}\ When interpreted in the context of our auxetic design principles, this result says that if a finite linkage $L$ in $R^d$ satisfies the strict auxetic prerequisites for $d$ pairs of vertices, then any affine
transform of $L$ will satisfy those prerequisites for the corresponding $d$ pairs. Nevertheless,
the intervals where an auxetic path for the associated periodic frameworks would be defined
may well differ.

\section{An infinite gallery}\label{sec: gallery}

Our gallery is dedicated to new examples in dimension three. We show that our
main auxetic design principle can be applied to an infinite series of finite linkages in $R^3$ of a rather elementary nature.  
We use minimally rigid linkages with one edge removed and adequate marking of three pairs of vertices. 

\begin{figure}[h]
\centering
 {\includegraphics[width=0.8\textwidth]{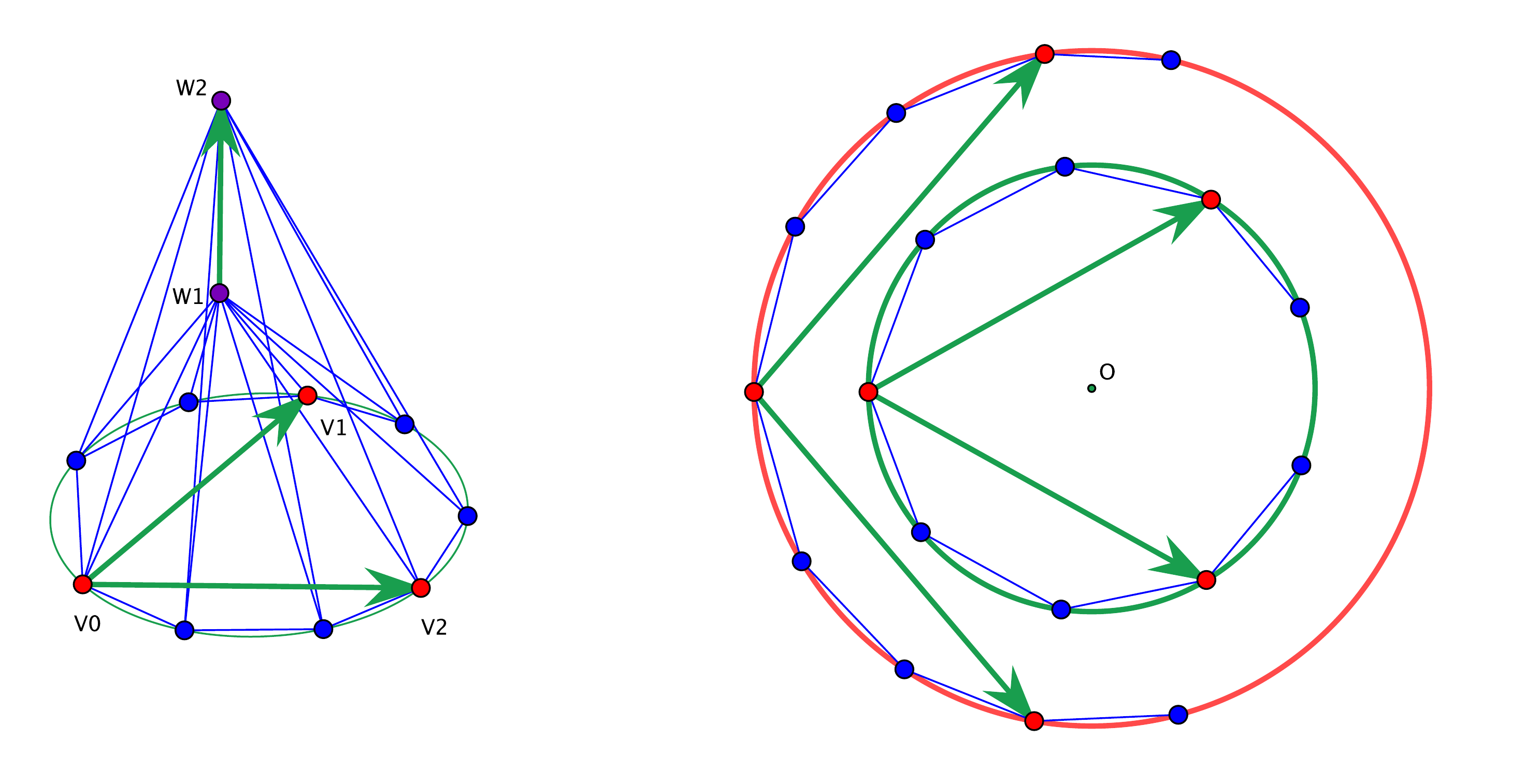}}
 \caption{The finite linkage $L_3$ with a diagram for the deformation effect in the horizontal plane.}
 \label{FigLk}
\end{figure}

\medskip \noindent
Let $k\geq 3$ be an integer and consider a regular polygon with $3k$ edges inscribed in the unit circle. The plane of the circle will be called horizontal.
We mark as $v_0,v_1,v_2$ three vertices which have from one to another exactly $n$ edges and
form an equilateral triangle. Then, we take two points $w_1,w_2$ on the vertical axis through
the center of the circle. The linkage denoted $M_k$ has $n=3k+2$ vertices and $m=9k$
edges and is obtained from the regular polygon by connecting $w_1$ and $w_2$ to all its vertices. When $w_1$ and $w_2$ are on opposite sides of the horizontal plane, we have the
one-skeleton of a convex polyhedron with triangular faces. By Cauchy's theorem, the linkage
is minimally rigid. By reflecting one part of the structure to the other side of the horizontal plane,
we see that minimal rigidity holds as well, hence $M_k$ is minimally rigid either way. We'll need the case with $w_1$ and $w_2$ on the same side, say above the horizontal plane, and we assume $w_1$ closer to the plane. The linkage to be used for obtaining a periodic framework is denoted
$L_k$ and is derived from $M_k$ by deleting an edge from the polygonal chain between $v_1$
and $v_2$ and marking the three pairs of vertices $(w_1,w_2)$, $(v_0,v_1)$ and $(v_0,v_2)$.
Figure~\ref{FigLk} shows $L_3$ next to the planar diagram used for describing its motion.

\noindent
The deformation mechanism of $L_k$ is easily represented since the vertices of the horizontal
polygon must remain coplanar: they are equally distanced from $w_1$, respectively $w_2$, hence
on the circle of intersection of two spheres. When following the deformation trajectory in the
direction of an augmenting radius for this intersection circle, the vector from $w_1$ to $w_2$
remains orthogonal to the circle plane and locally increases in length. Thus, for strict auxeticity,
we have to examine only the variation of the $2\times 2$ Gram matrix of the two vectors, say $\lambda_1$, $\lambda_2$, corresponding to the pairs $(v_0,v_1)$ and $(v_0,v_2)$. This is an elementary computation. With the radius $r$ of the circle as parameter, we find:

\begin{equation}\label{var}
\lambda_{1,2}=r(1-\cos{\theta},\pm \sin{\theta})
\end{equation}

\noindent 
where the expression of $\theta$ as a function of $r$ is given by

\begin{equation}\label{angle}
\theta=\theta(r)=2k \arcsin{(\frac{1}{r}\sin{\frac{\pi}{3k}})}
\end{equation}

\medskip \noindent
For the Gram matrix $\omega(r)=(\langle \lambda_i, \lambda_j \rangle)$, $1\leq i,j\leq 2$, it follows that $\frac{d\omega}{dr}(1)$ is positive definite. We conclude (via Corollary~\ref{strict})
that the one degree of freedom periodic framework associated to $L_k$ is a strictly auxetic periodic mechanism in a neighborhood of the initial position. Thus, with $k\geq 3$, we obtain an
infinite series of distinct auxetic periodic structures.

\medskip \noindent
{\bf Remark.} \ Considering that a small enough change in the placement of the vertices 
will maintain strict auxeticity,  the self-crossing resulting in the periodic framework from the planarity of the polygonal vertices in $L_k$ can be avoided by starting with a slightly distorted
version of $L_k$.
 
\section{Conclusion}
\label{sec: conclusion}

 The mathematical design principles presented here are based on the geometric foundations of
periodic auxetics introduced in \cite{borcea:streinu:geometricAuxetics:RSPA:arxiv:2015}. 
We have shown that, in spite of a rather conspicuous sparsity of auxetic designs in the existing literature, there are unlimited possibilities for generating periodic frameworks with auxetic capabilities in any dimension. The main procedure discussed in this paper converts any finite bar-and-joint mechanism in $R^d$ with adequate motion for $d$ marked pairs of joints into an auxetic periodic framework mechanism.  

\medskip
This work has brought to higher visibility a number of topics and features which deserve further
investigation. We propose them as open problems.

\medskip \noindent
(1)\ Identify the class of periodic graphs obtained by the passage from finite to periodic
described in Section~\ref{sec: F to P}.

\medskip \noindent
(2)\ Formulate precise and useful criteria for equivalence/non-equivalence of auxetic periodic mechanisms.

\medskip \noindent
(3)\ Determine and control the interval where a deformation path remains auxetic.

\paragraph{Acknowledgement.}The first author acknowledges partial support through NSF award no. 1319389 and the second author acknowledges partial support through NSF award no. 1319366. Both authors are partially supported through NIH Grant 1R01GM109456.  All statements, findings or conclusions contained in this paper are those of the authors and do not necessarily reflect the position or policy of the US Government. No official endorsement should be inferred.

\end{document}